\documentclass[12pt]{article}
\usepackage{amsmath,amsthm,amsfonts,latexsym,amssymb,enumerate,color}
\usepackage{graphics}
\usepackage{graphicx}
\usepackage{epsfig}

\newcommand\CC{{\mathbb C}}
\newcommand\GA{\mathcal A}

\newcommand\RR{{\mathbb R}}

\newcommand\GX{{\mathcal X}}

\def\beq{\begin{equation}}
\def\eeq{\end{equation}}
\newtheorem{thm}{Theorem}[section]

\newtheorem{lem}[thm]{Lemma}
\newtheorem{cor}[thm]{Corollary}
\newtheorem{rem}[thm]{Remark}

\newtheorem{ex}[thm]{Example}
\newcommand\beginpf{\noindent {\bf Proof:} \quad}

\newcommand\re{\mathop{\rm Re}\nolimits}
\newcommand\im{\mathop{\rm Im}\nolimits}

\newcommand\supp{\mathop{\rm supp}\nolimits}

\newcommand\rad{\mathop{\rm Rad}\nolimits}

\def\beginpf{\begin{proof}}
\def\endpf{\end{proof}}

\def\C{{\mathbb C}}

\newcommand\inside{\mathop{\rm int}\nolimits}
\renewcommand\phi{\varphi}
\newcommand\FB{\mathcal{F}\mathcal{B}}

\title{Estimates near the origin for functional calculus on analytic semigroups}
\author{ I. Chalendar\thanks{Universit\'e Paris Est Marne-la-Vall\'ee, 5 bd Descartes, Champs-sur-Marne,
77454 Marne-la-Vall\'ee, cedex 2, France  \protect\linebreak[3]
{\tt isabelle.chalendar@u-pem.fr}.},
\ J. Esterle\thanks{IMB, UMR 5251,  Universit\'e de Bordeaux, 351
cours de la Lib\'eration,
33405 Talence Cedex, France.
  \protect\linebreak[3]
{\tt jean.esterle@math.u-bordeaux1.fr}.}
\ and   J.R. Partington\thanks{School of
Mathematics,
University of Leeds, Leeds LS2 9JT, U.K.
\protect\linebreak[3]
{\tt J.R.Partington@leeds.ac.uk}.}
}

\begin{document}

\baselineskip18pt

\maketitle

\bibliographystyle{plain}

\begin{abstract}
This paper provides sharp lower estimates near the origin for the functional calculus $F(-uA)$ of a generator $A$ of 
an operator semigroup defined
on a sector;
here $F$ is given  as the Fourier--Borel transform of an analytic functional.
The results are linked to the existence of an identity element 
in the Banach algebra generated by the semigroup.
Both the quasinilpotent and non-quasinilpotent
cases are considered, and sharp results are proved extending many in the literature. 
\end{abstract}



\section{Introduction}

The purpose of this paper is to prove results concerning norm estimates in analytic
semigroups which complement the results proved in \cite{bottle} for
semigroups defined on $\RR^+$.
For good references on analytic semigroups we recommend the books \cite{EN,sinclair}:
we note that analytic semigroups of operators are norm-continuous in the open sector of the plane on which they are defined,
and thus also act as semigroups by multiplication on the Banach algebra that they generate.
Thus we may easily pass from the language of operators to the language of Banach algebras.

In \cite{BCEP} the following result was proved for semigroups defined on the right-hand half-plane $\CC_+$.
Here, $\rho$ denotes the spectral radius of an operator, and $\rad$ denotes the radical of an algebra.
  \begin{thm}\label{thm:zohra}
 Let $(T(t))_{t\in \CC_+}$ be an analytic non-quasinilpotent semigroup in a Banach algebra. Let ${\cal A}_T$ be the closed subalgebra  generated by $(T(t))_{t\in \CC_+}$ and let $\gamma >0$.
   If there exists $t_0>0$ such that 
   $$\sup_{t \in \CC_+, \vert t \vert \le t_0}\rho(T(t)-T((\gamma+1)t))<2$$
   then ${\cal A}_T/\rad {\cal A}_T$ is unital, and the generator of $(\pi(T(t))_{t>0}$ is bounded, where $\pi:{\cal A}_T\to {\cal A}_T/\rad {\cal A}_T$ denotes the canonical surjection.
   \end{thm}

This can be seen as a lower estimate for a functional calculus in $\GA_T$, determined by
$F(-A)=T(t)-T((\gamma+1)t)$, where
 $F: s \to e^{-s}-e^{-(\gamma+1)s}$ is the Laplace transform of the atomic measure $\delta_1-\delta_{\gamma+1}$.

This approach was taken in \cite{bottle} for semigroups defined on $\RR_+$, and very general results were proved
for both the quasinilpotent and non-quasinilpotent cases, providing extensions of results in \cite{EM,hille,kalton} and elsewhere.
For a detailed history of the subject, we refer to \cite{batty}.

Fewer results are available for analytic semigroups, 
and virtually nothing involving a general functional calculus,
although some dichotomy results are given in \cite{batty}.
To prove  more general results for analytic semigroups requires the notion of the Fourier--Borel transform of
a distribution acting on analytic functions, and in Section~\ref{sec:2} we define these transforms
and the associated functional calculus.

In Section~\ref{sec:3} we derive results in the non-quasinilpotent case, using the properties of the characters
defined on the algebra $\GA_T$, putting Theorem~\ref{thm:zohra} in a much more general
context.
Note that the results we prove are sharp, as is shown in Example~\ref{ex:sharp}.

Finally, the more difficult case of quasinilpotent semigroups is treated in Section~\ref{sec:4},
adapting the complex variable methods introduced in \cite{bottle}. Here there are additional
technical difficulties involved in defining the functional calculus, since we now work with measures supported on 
compact subsets of $\CC$.

\section{Analytic semigroups and functional calculus}
\label{sec:2}

For $0<\alpha< \pi/2$ let $S_\alpha$ denote the sector
\[
S_\alpha: = \{ z \in \CC_+: |\arg z| < \alpha\}.
\]
Let $H(S_\alpha)$ denote the Fr\'echet space of analytic functions on $S_\alpha$, endowed with the topology of local uniform convergence; thus, if $(K_n)_{n \ge 1}$ is an increasing sequence of
compact subsets of $S_\alpha$ with $\bigcup_{n \ge 1} K_n=S_\alpha$, we may specify the topology by the seminorms
\[
\|f\|_n:= \sup \{|f(z)|: z  \in K_n\}.
\]
Now let $\phi: H(S_\alpha) \to \CC$ be a continuous linear functional, in the sense that there is an index $n$ and a constant $M>0$ such that 
$|\langle f,\phi \rangle| \le M \|f\|_n$ for all $f \in H(S_\alpha)$. These
are sometimes known as {\em analytic functionals} \cite{hormander}.

We define the {\em Fourier--Borel transform\/} of $\phi$ by
\[
\FB(\phi)(z) = \langle e_{-z}, \phi \rangle,
\]
for $z \in \CC$, where $e_{-z}(\xi)=e^{-z\xi}$ for $\xi \in S_\alpha$.
This is an analogue of the Laplace transform, and is given under that name in \cite[Sec. 4.5]{hormander}.
We follow the terminology of \cite{FB}.

If $\phi \in H(S_\alpha)'$, as above, then by the Hahn--Banach theorem, it can be extended to a functional on $C(K_n)$, which we still write as $\phi$, and is thus
given by a Borel measure $\mu$ supported on $K_n$.

That is, we have
\[
\langle f,\phi\rangle = \int_{S_\alpha} f(\xi) \, d\mu(\xi),
\]
where $\mu$ (which is not unique) is a compactly supported measure. 
For example, if $\langle f,\phi \rangle = f'(1)$, then
\beq\label{eq:derivdist}
\langle f,\phi \rangle = \frac{1}{2\pi i} \int_C \frac{f(z) \, dz}{(z-1)^2},
\eeq
where $C$ is any sufficiently small circle surrounding the point 1.

Note that
 \[
\FB(\phi)(z)=\int_{K_n} e^{-z\xi} d\mu(\xi),
\]
and thus it is an entire function of $z$ satisfying 
 $ \sup_{\re z > r} |\FB(\phi)(z)| \to 0$ as $ r \to \infty$.
We shall sometimes find it convenient to
use the alternative notation $\FB(\mu)$.

Now let $T:=(T(t))_{t \in S_\alpha}$ be an analytic semigroup  on a  Banach space $\GX$, with infinitesimal generator $A$. Let $\phi\in H(S_\alpha)'$ and let $F=\FB(\phi)$.

We may thus define, formally to start with,
\[
F(-A)=\langle T,\phi \rangle = \int_{S_\alpha} T(\xi) \, d\mu(\xi),
\]
which is well-defined as a Bochner integral in $\GA_T$.
It is easy to verify that the definition is independent of the choice of $\mu$ representing $\phi$.

Moreover, if $u \in S_{\alpha-\beta}$, where $\supp \mu \subset  { S_\beta}$ and $0<\beta<\alpha$, then
we may also define 
\[
F(-uA) = \int_{S_\beta} T(u\xi) \, d\mu(\xi),
\]
since $u\xi$ lies in $S_\alpha$.

 In the following, a symmetric measure is a measure such that
$\mu(\overline S)=\overline{\mu(S)}$ for $S \subset S_\alpha$.  A symmetric measure will have a Fourier--Borel transform $F$ 
satisfying $F(z)=\widetilde F(z):=\overline{F(\overline z)}$ for all $z \in \CC$.

\section{The non-quasinilpotent case}
\label{sec:3}

For a semigroup $(T(t))_{t \in S_\alpha}$ with generator $A$, we write $\GA_T$ for the
commutative Banach algebra generated by the elements of the semigroup, and 
$\widehat \GA_T$ for its character space (Gelfand space).

The following result 
 is proved for semigroups on $\RR_+$ in
\cite[Lem.~3.1]{bottle} (see also
 \cite[Lem. 3.1]{EM} and \cite[Lem. 3.1]{BCEP}). It enables us to regard $A$ itself
as an element of $C(\widehat \GA_T)$ by defining an appropriate value $\chi(A)=-a_\chi$ for each $\chi \in \widehat \GA_T$.

\begin{lem}\label{lem:em31}
For a strongly continuous and eventually norm-continuous semigroup $(T(t))_{t >0}$ and a nontrivial character $\chi \in \widehat \GA_T$ there is a unique $a_\chi \in \CC$
such that $\chi(T(t))=e^{-ta_\chi}$ for all $t>0$. Moreover, the mapping
$\chi \mapsto a_\chi$ is continuous, and $\chi(F(-uA))=F(ua_\chi)$ in the case that $F$
is the Laplace transform of a measure $\mu$ on $\RR_+$.
\end{lem}

A similar result holds for analytic semigroups, with the same proof, where now $F$ is the Fourier--Borel transform of
a distribution (defined above).

The following theorem extends \cite[Thm. 3.6]{BCEP}

\begin{thm}\label{thm:moregeneral}
Let  $0 < \beta < \alpha < \pi/2$.
Let $\phi \in H(S_\alpha)'$, induced by a symmetric measure $\mu \in M_c( { S_\beta})$ such that $\int_{ { S_\beta}} d\mu(z)=0$,
and let $F=\FB( \phi)$.  
Let $(T(t))_{t\in S_\alpha}=(\exp(tA))_{t \in S_\alpha}$ be an analytic non-quasinilpotent semigroup and let  
 $\GA_T$ be the subalgebra generated by
$(T(t))_{t\in S_\alpha}$.
If there exists $t_0>0$ such that 
\begin{equation}\label{eq:refwantsit}
 \sup_{t\in S_{\alpha-\beta},|t|\le t_0}\rho (F(-tA))< \sup_{z\in S_{\alpha-\beta}}|F(z)|,
 \end{equation}
then ${\GA_T}/\rad{\GA_T}$ is unital
and the generator of $\pi(T(t))_{t\in S_{\alpha}}$ is bounded, where $\pi: {\GA_T}\to {{\GA_T}/\rad({\GA_T})}$ denotes the canonical surjection.
\end{thm}

\beginpf
By the maximum principle, for each $\theta\in (-\pi/2,\pi/2)$,
$F$ attains its maximum absolute value
$M_\theta$, say, on the ray $R_\theta=\{z:\arg z = \theta\}$ and $M_\theta$ is an increasing function of $\theta$ on $[0,\pi/2)$. 
Moreover,
$M_\theta=M_{-\theta}$, since $\mu$ is symmetric.

Clearly there exists a $d>0$ such that 
 the maximum value of $F$ on each ray $R_\theta$ is attained at a point  $z$ such that
 $\re z \le d$.

By Lemma \ref{lem:em31}, the hypotheses of the theorem, including non-quasinilpotency,  imply that
for each $\chi\in \widehat\GA_T$ there
exists ${a_\chi}\in \C$ such that $\chi(T(t))=e^{-{a_\chi}t}$ for all
$t\in S_\alpha$, and hence $\chi(F(-tA))=F({a_\chi}t)$.
Moreover, we know from \eqref{eq:refwantsit} that
\[
|F({a_\chi}t)| <  \sup_{z\in S_{\alpha-\beta}}|F(z)|
\]
for all $t \in S_{\alpha-\beta}$ with $|t| \le t_0$.

If for any point $t$ in the sector $\{t\in S_{\alpha-\beta}: \, |t|\le t_0\}$ we have $\re {a_\chi} t > d$, and
$|\arg {a_\chi} t| \ge \alpha-\beta$, then 
\[
|F(\lambda {a_\chi}t)| \ge  \sup_{z\in S_{\alpha-\beta}}|F(z)|
\]
for some real $\lambda $ between $0$ and $1$, giving a contradiction. In particular, if $|\arg {a_\chi}t_0| \ge \alpha-\beta$, then
$\re {a_\chi}t_0 \le d$.

Now suppose that $0 \le \gamma=\arg {a_\chi}t_0 < \alpha-\beta$ (the other case is similar); then we know that
$\re( {a_\chi}t_0 e^{i(\alpha-\beta-\gamma)}) \le d$; writing ${a_\chi}t_0=re^{i\gamma}$ and
${a_\chi}t_0 e^{i(\alpha-\beta-\gamma)} = re^{i(\alpha-\beta)}$, we deduce that
$\re {a_\chi}t_0 \le d\cos\gamma / \cos (\alpha-\beta)$.
Hence, in all cases, we obtain
$\re {a_\chi}t_0 \le d / \cos (\alpha-\beta)$ and
\[
|\chi(T(t_0))| \ge \exp(-d/\cos(\alpha-\beta)),
\]
 and this holds for every $\chi\in \widehat{\GA}_T$. 
Hence $\widehat{\GA}_T$ is compact. Now by \cite[Thm. 3.6.3, 3.6.6]{rickart}, we conclude that
$\GA_T/\rad \GA_T$ is unital. Since the algebra generated by the norm-continuous semigroup $(\pi(T(t))_{t>0}$ is dense in ${\GA_T}/\rad {\GA_T},$ the generator of this semigroup is bounded.

\endpf

By considering the convolution of a functional $\phi \in H(S_\alpha)'$, possessing the Fourier--Borel transform $F$, and 
the functional $\widetilde \phi$, possessing the Fourier--Borel transform $\widetilde F$ (defined by the
convolution of defining measures), we 
obtain a functional whose transform is the product $F\widetilde F$, and may deduce the following result.

\begin{cor}
Let  $0 < \beta< \alpha < \pi/2$.
Let $\phi \in H(S_\alpha)'$, induced by a  measure $\mu \in M_c(S_\beta)$ such that $\int_{S_\beta} d\mu(z)=0$,
and let $F=\FB( \phi)$.
Let $(T(t))_{t\in S_\alpha}=(\exp(tA))_{t \in S_\alpha}$ be an analytic non-quasinilpotent semigroup  and let  
 $\GA_T$ be the subalgebra generated by
$(T(t))_{t\in S_\alpha}$.
If there exists $t_0>0$ such that 
\[ \sup_{t\in S_{\alpha-\beta},|t|\le t_0}\rho (F(-tA)\widetilde F(- t A))< \sup_{z\in S_{\alpha-\beta}}|F(z)||\widetilde F(z)|,\]
then ${\GA_T}/\rad{\GA_T}$ is unital
and the generator of $\pi(T(t))_{t\in S_{\alpha}}$ is bounded, where $\pi: {\GA_T}\to {{\GA_T}/\rad({\GA_T})}$ denotes the canonical surjection.
\end{cor}

A similar proof gives the following result.

\begin{thm}
Let  $0 < \alpha < \pi/2$.
Let $\phi \in H(S_\alpha)'$, induced by a  symmetric measure $\mu \in M_c(S_\alpha)$ such that $\int_{S_\alpha} d\mu(z)=0$,
and let $F=\FB( \phi)$. 
Let $(T(t))_{t\in S_\alpha}=(\exp(tA))_{t \in S_\alpha}$ be an analytic non-quasinilpotent semigroup  and let  
 $\GA_T$ be the subalgebra generated by
$(T(t))_{t\in S_\alpha}$.
If there exists $t_0>0$ such that 
\[ \rho (F(-tA))< \sup_{x>0}|F(x)|,\]
for all $0 < t \le t_0$,
then ${\GA_T}/\rad{\GA_T}$ is unital
and the generator of $\pi(T(t))_{t\in S_{\alpha}}$ is bounded, where $\pi: {\GA_T}\to {{\GA_T}/\rad({\GA_T})}$ denotes the canonical surjection.
\end{thm}

\beginpf
The proof is similar to that of Theorem~\ref{thm:moregeneral}, and we adopt
the same notation. We cannot have $\re a_\chi t_0 > d$ for any $\chi \in \widehat{\GA}_T$,
as then there would be a $\lambda \in (0,1)$ with 
\[
|\chi(F(-\lambda t_0 A))|=|F(\lambda a_\chi t_0)| \ge \sup_{x>0} |F(x)|.
\]
 We conclude
that $|\chi(T(t_0))| \ge \exp(-d)$ for all $\chi \in \widehat{\GA}_T$. 
The proof is now concluded as before.

\endpf

Again there is an immediate corollary for general $\phi$.

\begin{cor}
Let  $0 < \alpha < \pi/2$.
Let $\phi \in H(S_\alpha)'$, induced by a  measure $\mu \in M_c(S_\alpha)$ such that $\int_{S_\alpha} d\mu(z)=0$,
and let $F=\FB( \phi)$. 
Let $(T(t))_{t\in S_\alpha}=(\exp(tA))_{t \in S_\alpha}$ be an analytic non-quasinilpotent semigroup  and let  
 $\GA_T$ be the subalgebra generated by
$(T(t))_{t\in S_\alpha}$.
If there exists $t_0>0$ such that 
\[ \rho (F(-tA)\widetilde F(-tA))< \sup_{x>0}|F(x)|^2,\]
for all $0 < t \le t_0$,
then ${\GA_T}/\rad{\GA_T}$ is unital
and the generator of $\pi(T(t))_{t\in S_{\alpha}}$ is bounded, where $\pi: {\GA_T}\to {{\GA_T}/\rad({\GA_T})}$ denotes the canonical surjection.
\end{cor}

\begin{ex} \label{ex:sharp}
In the Banach algebra
$\GA=C_0[0,1]$ consider the semigroup $T(t):x \mapsto x^t$ for $t \in S_\alpha$. Clearly there is no convergence as $t$ approaches $0$.

For $x \in (0,1]$ (which can be identified with the Gelfand space of $\GA$) let $F=\FB(\mu)$ and
\[
F(-tA)(x) = \int_{S_\alpha} x^{t\xi} \, d\mu(\xi) = \int_{S_\alpha} e^{t\xi \log x} \, d\mu(\xi),
\]
where $\mu \in M_c(S_\alpha)$, supposing that  $\int_{S_\alpha} d\mu(z)=0$ and that $F$ is real for $x>0$.

Thus $F(-tA)(x)=F(-t \log x)$ and 
\[
\rho(F(-tA))=\|F(-tA)\|= \sup_{x \in (0,1]}|F(-t\log x)|=\sup_{r>0}|F(tr)|.
\]
Clearly  
\[ \sup_{t\in S_\alpha,|t|\le t_0}\rho (F(-tA))= \sup_{t\in S_\alpha}|F(z)|
\]
for all $t_0>0$. This shows that the hypothesis \eqref{eq:refwantsit} of Theorem~\ref{thm:moregeneral} is sharp.
\end{ex}

\begin{rem}
In
 \cite[Thm.~3.3]{bottle} it is shown that 
if $F$ is the Laplace transform of a real compactly-supported measure $\mu$ with $\int_0^\infty d\mu=0$, then
the condition 
\[
\rho(F(-u_k A)) < \sup_{t>0} |F(u_kt)|
\]
for a (real) sequence $u_k \to 0$ implies that the  algebra $ \GA_T$ possesses an exhaustive sequence of idempotents $(P_n)_{n \ge 1}$ (i.e., $P_n^2=P_nP_{n+1} = P_n$ for all $n$ and for every $\chi \in \widehat \GA_T$ there is a $p$ such that
$\chi(P_n)=1$ for all $n \ge p$),  
such that each semigroup $(P_n T(t))_{t>0}$ has a bounded generator.

The directly analogous result for analytic semigroups does not hold. For example, we may consider a modification of Example
\ref{ex:sharp}, by taking $\GA=C_0[0,1]$ and the semigroup $T(t): x \to x^t= \exp(t \log x)$ for $x \in (0,1]$ 
and $t \in S_{\pi/2}$ with generator $A$, say. Let
$F(z)=e^{-z}-e^{-2z}=\FB(\delta_1-\delta_2)$.
Now let $\omega=\exp(i\pi/6)$ and define
 $\widetilde T(t)=T(\omega t)$.
Thus $\widetilde T(t)$ is an analytic semigroup in $S_{\pi/3}$, with generator $\widetilde A=\omega A$.
Define $\widetilde F(z)=F(z/\omega)$, the Fourier--Borel transform of a measure supported in $S_{\pi/3}$.
Now, for $u>0$,  
\[
\rho(\widetilde F(-u\widetilde A))=\rho(F(-uA))=\rho(x^u-x^{2u})=1/4;
\]
however,
\[
\sup_{t >0} |\widetilde F(tu)|=\sup\{|F(tu/\omega)|: t > 0\}> 1/4,
\]
 by the maximum principle
applied to $F$ on the sector $S_{\pi/6}$ (numerically, the supremum is about 0.29).
\end{rem}

\section{The quasinilpotent case}
\label{sec:4}

We suppose throughout this section that $(T(t))_{t \in S_\alpha}$ is a  nontrivial strongly continuous semigroup of 
uniformly bounded
quasi\-nilpotent operators acting on a 
Banach space $(\GX,\|.\|)$.

\subsection{Preliminary results}

We shall need the following theorem from complex analysis, which was proved in \cite[Thm. 2.2]{bottle}.

\begin{thm}\label{thm:Jordan}
Let $f: \overline{\CC_+} \to \CC$ be a continuous bounded nonconstant function, analytic on $\CC_+$, such that
$f([0,\infty)) \subset \RR$, $f(0)=0$,   satisfying the condition $\lim_{ { x \to \infty},{x \in \RR}} f(x)=0$.

Suppose that $b > 0$ is such that $f(b) \ge |f(x)| $ for all $x \in [0,\infty)$.
Then there exist $a_1,a_2 \in \CC_+$, $a_0 \in (b,a_1)$ and $a_3 \in i\RR$ with $\im a_j > 0$ for  $j=1,2,3$, and $\im a_2=\im a_3$,
and a simple piecewise linear 
Jordan curve $\Gamma_1$ joining $a_1$ to $a_2$ in the upper right half-plane
$\{z \in \CC: \re z>0, \im z >0\}$ and $\delta > 0$ such that

(i) $|f(z)| \ge f(b) + \delta |z-b|^m$ for all $z \in [b,a_1]$, where $m$ (even) is the smallest
positive integer with $f^{(m)}(b) \ne 0$;

(ii) $|f(z)| > |f(a_0)|$ for all $z \in \Gamma_1 \cup [a_2,a_3]$.
\end{thm}

\begin{figure}[htbp]
  \begin{center}
\begin{picture}(0,0)%
\includegraphics[scale=0.5]{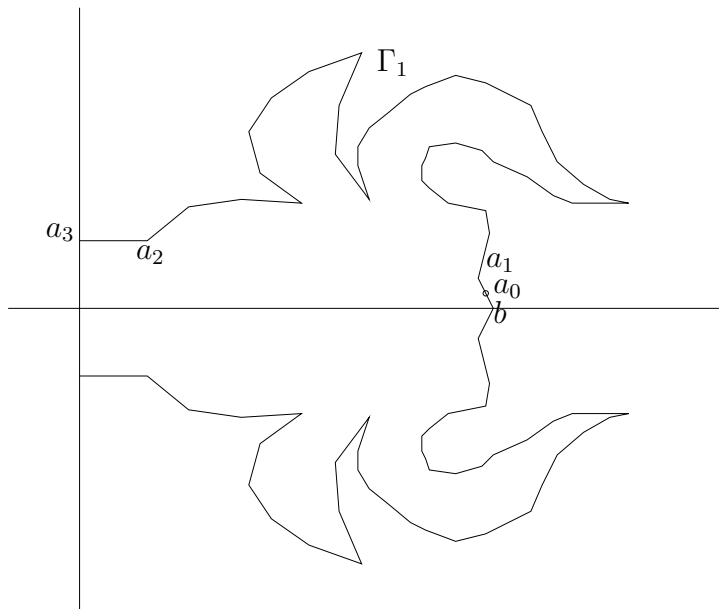}%
\end{picture}%
\setlength{\unitlength}{2072sp}%
\begingroup\makeatletter\ifx\SetFigFont\undefined%
\gdef\SetFigFont#1#2#3#4#5{%
  \reset@font\fontsize{#1}{#2pt}%
  \fontfamily{#3}\fontseries{#4}\fontshape{#5}%
  \selectfont}%
\fi\endgroup%
\begin{picture}(8529,7224)(1879,-6823)
\put(3421,-2626){\makebox(0,0)[lb]{\smash{{\SetFigFont{12}{14.4}{\rmdefault}{\mddefault}{\updefault}{\color[rgb]{0,0,0}$a_2$}%
}}}}
\put(2341,-2356){\makebox(0,0)[lb]{\smash{{\SetFigFont{12}{14.4}{\rmdefault}{\mddefault}{\updefault}{\color[rgb]{0,0,0}$a_3$}%
}}}}
\put(6301,-376){\makebox(0,0)[lb]{\smash{{\SetFigFont{12}{14.4}{\rmdefault}{\mddefault}{\updefault}{\color[rgb]{0,0,0}$\Gamma_1$}%
}}}}
\put(7606,-2716){\makebox(0,0)[lb]{\smash{{\SetFigFont{12}{14.4}{\rmdefault}{\mddefault}{\updefault}{\color[rgb]{0,0,0}$a_1$}%
}}}}
\put(7696,-3031){\makebox(0,0)[lb]{\smash{{\SetFigFont{12}{14.4}{\rmdefault}{\mddefault}{\updefault}{\color[rgb]{0,0,0}$a_0$}%
}}}}
\put(7696,-3391){\makebox(0,0)[lb]{\smash{{\SetFigFont{12}{14.4}{\rmdefault}{\mddefault}{\updefault}{\color[rgb]{0,0,0}$b$}%
}}}}
\end{picture}%
   \caption{The curve given in Theorem \ref{thm:Jordan}}
    \label{fig:g1}
  \end{center}
\end{figure}

Moreover, the following formula is easily derived.

\begin{lem}\label{lem:formule}
Let $(T(t)=(e^{tA}))_{t\in S_{\alpha}}$ be an analytic   semigroup on $S_{\alpha},$ and let $\phi \in { H}(S_{\alpha})'$ be  a   functional represented by a measure $\mu$ supported on a compact set  $K \subset S_{\alpha}$,
so that we have
\[
\langle f , \phi \rangle  =\int_{K}f(\zeta)d\mu(\zeta) \qquad (f \in {\mathcal H}(S_{\alpha})).
\]
Assume that $\sup_{t \in S_{\beta}\cap D(0,1)}\Vert T(t)\Vert <+\infty$ for every $\beta \in [0, \alpha),$ so that the Bochner integral 
\[
 \int_K T(s\zeta u)e^{(s-1)\zeta \lambda}\zeta \, d\mu(\zeta)
\]
 is well-defined in ${\mathcal B}(\GX)$ for every $\zeta \in K$, $\lambda \in \CC$, and $u\in\CC$ with $uK \subset S_{\alpha}$.
Then, for $F=\FB(\phi)$ we have
$$F(-uA) - F(\lambda)I= (uA+\lambda I)\left ( \int_{0}^1\left [ \int_K T(s\zeta u)e^{(s-1)\zeta \lambda}\zeta d\mu(\zeta)\right ] ds\right ).$$
Also, if $x \in D(A),$ we have

$$F(-uA)x - F(\lambda)x= \left ( \int_{0}^1\left [ \int_K T(s\zeta u)e^{(s-1)\zeta \lambda}\zeta d\mu(\zeta)\right ] ds\right )(uA+\lambda I)x.$$
\end{lem}

\beginpf
Using Fubini's theorem and the  fundamental theorem of the calculus, we have
\begin{eqnarray*}
\int_{0}^1\left [ \int_K T(s\zeta u)e^{(s-1)\zeta \lambda}\zeta d\mu(\zeta)\right ] ds
&=& \int_K (uA+\lambda I)^{-1} (T(\zeta u)-e^{-\zeta\lambda}I) \, d\mu(\zeta)\\
&=& (uA+\lambda I)^{-1} (F(-uA)-F(\lambda)I),
\end{eqnarray*}
from which the result is clear.
\endpf

\begin{cor}\label{cor:ruby}
Let  $(T(t)=(e^{tA}))_{t\in S_{\alpha}}$ be a bounded analytic   semigroup on $S_{\alpha}$.
Then there exists a constant $C>0$ such that
\[
\|(F(-uA)-F(\lambda) I)(uA+\lambda I)^{-1}\|
\le
C \int_K |\zeta| \, d|\mu| (\zeta),
\]
for all $u$ with
 $uK \subset S_\alpha$ and $|\lambda| \le 1$.
\end{cor}

\beginpf
This follows immediately from 
Lemma~\ref{lem:formule} since $|e^{(s-1)\zeta\lambda}| $ is uniformly bounded for $\zeta\in K$ and 
$s \in [0,1]$.

\endpf

We shall require one further preliminary result, which is well-known for $z$ in the closed left half-plane,
and was proved in \cite{bottle}.

\begin{lem}\label{lem:resolvbound}
Let $(T(t))_{t \in S_\alpha}$ be bounded analytic quasinilpotent
semigroup with generator $A$, and let $0<\beta<\alpha$. Then
$\|(A+zI)^{-1}\|$ is uniformly bounded on the set
\[
\left \{z \in \CC: \frac{\pi}{2} - (\alpha-\beta) \le |\arg z| \le \pi \right\}.
\]
\end{lem}

\beginpf
It is standard that  $\|z(A+zI)^{-1}\|$ is bounded
in the given set for any bounded analytic semigroup. Moreover, quasinilpotency implies that $A$ is invertible and then
$\|(A+zI)^{-1}\|$ is bounded near 0.
\endpf

\subsection{Estimates for quasinilpotent semigroups}


The following theorem is a version of \cite[Thm.~2.5]{bottle} for analytic semigroups. The proof is
based on the same   ideas from complex analysis with several significant modifications.
For $0<\beta< \pi/2$ let
\[
V_\beta:= \{ z \in \CC_+: |\arg z| \le \beta\}.
\]

\begin{thm}\label{thm:leffe}
Let $(T(t))_{t \in S_\alpha}$ be a nontrivial  bounded semigroup of quasinilpotent operators and let $F$ be the Fourier--Borel transform of a 
symmetric measure $\mu$ supported on a compact set $K \subset V_\beta$ 
for some  $\beta$ with $0< \beta < \alpha$,
such that $\int_K \, d\mu=0$.

Then 
there is an $\eta >0 $ such that
\[
\|F(-uA)\| > \sup_{t > 0} |F(t)| 
\]
for all $u \in S_{\alpha-\beta}$ with $|u| \le \eta$.
\end{thm}


\beginpf																																																																																																															
Note that $F$ is real on the positive real axis, since $\mu$ is symmetric.
 Let $b>0$ be such that
$|F(x)| \le |F(b)|$ for all $x \ge 0$. By considering $-\mu$ instead of $\mu$, if necessary, we may
suppose that $F(b) > 0$.

By Corollary~\ref{cor:ruby}
\[
\left \| F(-uA)(uA+\lambda I)^{-1} \right \| \ge \left\|F(\lambda)(uA+\lambda I)^{-1} \right\| -C \int_K |\zeta| \, d|\mu| (\zeta).
\]
for $u \in S_{\alpha-\beta}$ and $|\lambda| \le 1$, and hence
\[
\|F(-uA)\| \ge |F(\lambda)|-\frac{C}{\left\|(uA+\lambda I)^{-1}\right\|}\int_K |\zeta| \, d|\mu| (\zeta) .
\]

Suppose that there exists a $u \in S_{\alpha-\beta}$ with $|u|<1$ such that
$\|F(-uA)\| \le F(b)$, and consider the simple Jordan curve
\[
\Gamma:=[b,a_1] \cup \Gamma_1 \cup [a_2,a_3] \cup [a_3,\overline{a_3}] \cup [\overline{a_3},\overline{a_2}]
\cup \overline{\Gamma_1} \cup [\overline a_1, b],
\]
where $\Gamma_1,a_1,a_2,a_3$ are defined as in Theorem~\ref{thm:Jordan}, taking $f=F$ (see Figure
\ref{fig:g1}).

We now make various estimates of $\left\|(uA+\lambda I)^{-1}\right\|$ for $\lambda$ on three different parts of $\Gamma$.

1) By (i) of Theorem \ref{thm:Jordan}, for $\lambda \in [b,a_1] \cup [\overline{a_1},b]$ we have
\begin{eqnarray*}
F(b) \ge \|F(-uA)\| &\ge& |F(\lambda)| - \frac{C}{\left\|(uA+\lambda I)^{-1}\right\|}\int_K |\zeta| \, d|\mu| (\zeta) \\
&\ge& F(b) + \delta|\lambda-b|^m - \frac{C}{\left\|(uA+\lambda I)^{-1}\right\|}\int_K |\zeta| \, d|\mu| (\zeta).
\end{eqnarray*}
Hence we obtain
\beq\label{eq:est1}
\left\|(uA+\lambda I)^{-1}\right\| \le \frac{C}{\delta|\lambda-b|^m}\int_K |\zeta| \, d|\mu| (\zeta).
\eeq

2) By (ii) of Theorem \ref{thm:Jordan}, for $\lambda \in \Gamma_1 \cup [a_2,a_3] \cup [\overline{a_3},\overline{a_2}] \cup \overline{\Gamma_1}$ we have
\begin{eqnarray*}
F(b) \ge \|F(-uA)\| &\ge& |F(\lambda)| - \frac{C}{\left\|(uA+\lambda I)^{-1}\right\|}\int_K |\zeta| \, d|\mu| (\zeta)\\
&\ge& |F(a_0)| - \frac{C}{\left\|(uA+\lambda I)^{-1}\right\|}\int_K |\zeta| \, d|\mu| (\zeta).
\end{eqnarray*}
It follows that
\beq\label{eq:est2}
\left\|(uA+\lambda I)^{-1}\right\| \le \frac{C}{|F(a_0)|-F(b)}\int_K |\zeta| \, d|\mu| (\zeta).
\eeq


3)  By
Lemma~\ref{lem:resolvbound}  there is a constant $C'>0$ depending only on
the semigroup such that
\beq\label{eq:est3}
\left\| \left(A+\frac{\lambda}{u}I \right)^{-1}\right\| \le C' 
\eeq
for $\lambda \in [a_3,\overline{a_3}] \subset i\RR$.

We can now provide estimates for the quantity $\left\| (\lambda-b)^m \left( A+\frac{\lambda}{u}I \right)^{-1} \right\|$ 
for $\lambda$ on $\Gamma$. Let $R=\max_{\lambda \in \Gamma}|\lambda-b|$.

By (\ref{eq:est1})
\[
\left\| (\lambda-b)^m \left( A+\frac{\lambda}{u}I \right)^{-1} \right\| \le \frac{C|u|}{\delta} \int_K |\zeta| \, d|\mu| (\zeta)
\]
for all $\lambda \in [b,a_1] \cup [\overline{a_1},b]$.

By (\ref{eq:est2})
\[
\left\| (\lambda-b)^m \left( A+\frac{\lambda}{u}I \right)^{-1} \right\| \le \frac{C|u|R^m}{|F(a_0)|-F(b)}\int_K |\zeta| \, d|\mu| (\zeta)
\]
for all $\lambda \in \Gamma_1 \cup [a_2,a_3] \cup [\overline{a_3},\overline{a_2}] \cup \overline{\Gamma_1}$.

By (\ref{eq:est3})
\[
\left\| (\lambda-b)^m \left( A+\frac{\lambda}{u}I \right)^{-1} \right\| \le C'R^m 
\]
for all $\lambda \in [a_3,\overline{a_3}]$.

Since $0 < |u| \le 1$, for all $z \in \Gamma \cup \inside \Gamma$ we have
\[
\left\|  \left( A+\frac{z}{u}I \right)^{-1} \right\| \le \frac{M}{|z-b|^m},
\]
by the maximum modulus principle,
where
\[
M= \max \left(
\frac{C}{\delta} \int_K |\zeta| \, d|\mu| (\zeta), 
\frac{CR^m}{|F(a_0)|-F(b)}\int_K |\zeta| \, d|\mu| (\zeta),
C'R^m \right) .
\]

Since by hypothesis $F(0)=0$, there is an $r \in (0,b)$ such that
\[
\sup_{|z| \le r} |F(z)| < F(b).
\]
Taking $r$ sufficiently small, we have
$\overline{D(0,r)}\cap \Gamma \cap \CC_+ = \emptyset$, and then $\overline{D(0,r)} \cap \overline{\CC_+}
\subset \Gamma \cup \inside \Gamma$.

Now if $z \in \overline{D(0,r)}$ with $\re z >0$, we have
$|z-b| \ge b-r$, and thus we have
\beq\label{eq:estarm}
\left\| \left( A + \frac{z}{u}I \right)^{-1} \right\| \le \frac{M}{|z-b|^m} \le \frac{M}{(b-r)^m}.
\eeq
Also,  by Lemma~\ref{lem:resolvbound} there is an $M'>0$ such that
\[
\sup_{z \in E} \left\| \left( A + zI \right)^{-1} \right\| \le M',
\]
where 
\[
E=\left \{z \in \CC: \frac{\pi}{2} - (\alpha-\beta) \le |\arg z| \le \pi \right\}.
\]
Hence 
\begin{equation}\label{eq:zincminus}
\left\|\left (A+\frac{z}{v}I \right)^{-1} \right\| \le M'
\end{equation}
for $\re z \le 0$ and $v \in S_{\alpha-\beta}$.

Now, since by Liouville's theorem the function $\lambda \mapsto \left\| \left( A + \lambda I \right)^{-1} \right\| $ is unbounded on $\CC$,
and so there is a $\lambda \in \CC$ with 
\[
\left\|\left (A+ \lambda I \right)^{-1} \right\| > \max\left\{ M', \frac{M}{(b-r)^m}\right\}.
\]
It follows that there is an $\eta>0$ such that
  for all $u \in S_{\alpha-\beta}$ with $|u| \le \eta$ we have
\[
\left\|\left (A+\frac{z}{u}I \right)^{-1} \right\| > \max\left\{ M', \frac{M}{(b-r)^m}\right\}
\]
  for some $z=\lambda u \in \overline{D(0,r)}$. Note that $z \in   \CC_+$ by \eqref{eq:zincminus}.

This contradicts \eqref{eq:estarm}, and it follows that  
\[
\|F(-uA)\|> F(b) \qquad \hbox{for all} \quad u \in S_{\alpha-\beta}, \quad |u| \le \eta.
\]

\endpf

\begin{ex}
The analytic distribution $\phi: f \mapsto f'(1)$ discussed earlier (cf.  \eqref{eq:derivdist})
leads to $F(z)=-ze^{-z}$ and $F(-uA)=uAT(u)$.
Theorem \ref{thm:leffe} now asserts that $\|uAT(u)\| \ge 1/e$ for
sufficiently small
$u \in S_{\alpha-\beta}$, where $\beta$ corresponds to the support of $\phi$ and can 
therefore be taken arbitrarily small. 
This may be seen in the context of the result of Hille \cite{hille} (see also \cite[Thm. 10.3.6]{HP})
that if $\limsup_{t \to 0+} t \|AT(t)\| < 1/e$, then $A$
is bounded, and hence the semigroup cannot be quasinilpotent.
There is a  further discussion of this example in \cite[Thm. 3.12, Thm 3.13]{BCEP}.
\end{ex}

We have a corollary for measures that are not necessarily symmetric.
If $\mu$ is now a complex compactly-supported measure,
with Fourier--Borel transform $F$,  then we write $\widetilde F(z)=\overline{F(\overline z)}$; then $\widetilde F$
is also an entire function, indeed, the Fourier--Borel transform of $\overline{\mu}$.

\begin{cor}\label{cor:26z}
Let $(T(t))_{t \in S_\alpha}$ be a nontrivial  bounded analytic semigroup of quasinilpotent operators and let $F$ be the Fourier--Borel transform of a 
 measure $\mu$ supported on a compact set $K \subset V_\beta$ 
for some  $\beta$ with $0< \beta < \alpha$,
such that $\int_K \, d\mu=0$.

Then 
there is an $\eta >0 $ such that
\[
\|F(-uA)\widetilde F(-uA)\| \ge \sup_{t > 0} |F(t)|^2
\]
for all $u \in S_{\alpha-\beta}$ with $|u| \le \eta$.
\end{cor}

\beginpf
The result follows on applying
Theorem~\ref{thm:leffe} to the
 real measure $\nu:= \mu * \overline \mu$, whose Fourier--Borel transform satisfies
\[
\FB(\nu) (s) = F(s)\widetilde F(s).
\]
\endpf

\section*{Acknowledgements}

The authors thank the Institut Camille Jordan of  Universit\'e Lyon 1 and the IMB of the Universit\'e de
Bordeaux for their hospitality and financial support. They are also grateful to the referee for many useful
suggestions.

\end{document}